\documentclass[11pt]{amsart}
\usepackage[T1]{fontenc}
\usepackage{lmodern}
\usepackage{mathtools,amssymb,amsthm}
\usepackage{booktabs,array,flafter}
\usepackage{enumitem}
\usepackage{microtype,needspace}
\usepackage[a4paper,margin=30mm]{geometry}
\usepackage{hyperref}
\usepackage{float}
\usepackage{tikz}
\usetikzlibrary{arrows.meta,positioning,calc}
\hypersetup{
 pdftitle={The nonexistence of nine-point crescent configurations},
 pdfsubject={A computer-assisted proof using six-point cores and exact arithmetic},
 pdfkeywords={crescent configurations, distance multiplicities, distance geometry}
}

\newtheorem{theorem}{Theorem}[section]
\newtheorem{proposition}[theorem]{Proposition}
\newtheorem{lemma}[theorem]{Lemma}
\newtheorem{corollary}[theorem]{Corollary}
\newtheorem{computation}[theorem]{Finite computation}
\theoremstyle{definition}
\newtheorem{definition}[theorem]{Definition}
\newtheorem{procedure}[theorem]{Procedure}
\theoremstyle{remark}

\numberwithin{equation}{section}
\newcommand{\Q}{\mathbb Q}
\newcommand{\R}{\mathbb R}
\newcommand{\F}{\mathcal F}
\newcommand{\U}{\mathcal U}
\DeclareMathOperator{\rank}{rank}
\DeclareMathOperator{\spanQ}{span_{\Q}}
\setlist[enumerate]{label=\textup{(\roman*)},leftmargin=*,itemsep=3pt,topsep=5pt}

\title[The nonexistence of nine-point crescent configurations]{The nonexistence of nine-point crescent configurations}
\author{Yifei Wang}
\email{\textcolor{blue}{\href{mailto:yifeiwang.math@gmail.com}{yifeiwang.math@gmail.com}}}
\address{School of Mathematical Sciences and LPMC, Nankai University, Tianjin 300071, China}
\date{}

\begin{document}
\begin{abstract}
A planar crescent configuration consists of $n$ points, no three collinear
and no four concyclic, whose distinct distances have multiplicities
$1,2,\ldots,n-1$.  We prove by exact computation that no such configuration
exists for $n=9$.  The three pairs carrying the two least frequent distances
yield a six-point subset with at most six distances.  Rational linear
algebra on quartic polynomials separates the possible six-point distance
patterns into those with certified finite lists of normalized distance
tuples and a residual family.  Compatibility of overlapping subsets forces
at least one core into the first family, without classifying the residual
realizations.  The remaining three points
are then obtained by rational constructions over the field of the squared core
distances.  The multiplication matrices used to enumerate the core
distances have order at most twelve.  Exhaustive label assignments and exact
algebraic comparisons exclude all remaining extensions.
\end{abstract}
\maketitle

\section{Introduction}\label{sec:introduction}

Distance problems ask not only how many distinct lengths a point set
can determine, but also how often each length can occur.
Erd\H{o}s's distinct distances problem \cite{Erdos1946} concerns the
first question; Guth and Katz \cite{GuthKatz2015} proved that every
set of $n\ge2$ planar points determines at least $cn/\log n$ distinct
distances for an absolute constant $c>0$.  Prescribing multiplicities
asks a different question: which sequences can occur as the occurrence
counts of all distinct distances?  We study the sequence
$1,2,\ldots,n-1$, whose sum is exactly the number $\binom n2$ of
unordered pairs.

\begin{definition}\label{def:crescent}
	A finite set $P\subset\R^2$ is in \emph{general position} if no three
	of its points are collinear and no four are concyclic.
	It is a \emph{crescent configuration} if it is in general position
	and, writing $n=|P|$, it determines exactly $n-1$ distances, one of
	multiplicity $k$ for each $1\le k\le n-1$.
	Multiplicity counts unordered pairs of distinct points; the distances
	are not ordered by magnitude.
\end{definition}

The multiplicity problem appears in Erd\H{o}s's writings of the early
1980s \cite{Erdos1983,Erdos1986} and in his later account
\cite{Erdos1989}.  Without general position, equally spaced points on
a line give the required multiplicities for every $n$; equally spaced
points on a sufficiently short circular arc do so as well.
Erd\H{o}s conjectured that the general-position requirement changes
this behavior: there should be an integer $N$ such that no planar
crescent configuration exists for any $n\ge N$
\cite{Erdos1986,Erdos1989}.  The name \emph{crescent configuration}
was introduced by Burt, Goldstein, Manski, Miller, Palsson, and Suh
\cite{Burt2016}.

The small cases were resolved by a succession of constructions.
The cases $n=2,3,4$ are elementary: besides two distinct points and a
nonequilateral isosceles triangle, a suitably chosen isosceles triangle
together with its circumcenter gives four points with multiplicities
$1,2,3$.  Erd\H{o}s's 1983 account \cite{Erdos1983} describes
Pomerance's five-point example and records the existence of a six-point
example.  Seven-point constructions were obtained independently by
Andy Liu and Pal\'asti \cite{Erdos1986,Liu1986,Palasti1987}.
Pal\'asti subsequently constructed an eight-point example
\cite{Palasti1989}, establishing existence for every $2\le n\le8$.

The triangular lattice has played a central role in these
constructions and in the subsequent searches.
Pal\'asti's lattice examples \cite{Palasti1989} use the equilateral
triangular lattice, which, up to a similarity, is
\[
\left\{\left(a+\frac b2,\frac{\sqrt3}{2}b\right):
a,b\in\mathbb Z\right\}.
\]
A lattice displacement indexed by $(a,b)$ has squared length
$a^2+ab+b^2$, so equality of distances can be tested with integers.
Pal\'asti's eight-point configuration fits in a $37$-point patch
\cite{Fish2020}.  Burt et al.\ \cite{Burt2016} exhaustively searched
a $91$-point hexagonal region of the triangular lattice and found no
nine-point configuration there.  This is a finite-region exclusion,
not a nonexistence theorem for the whole lattice or the plane.
The general problem has no lattice hypothesis: a lattice search
restricts the available positions, whereas an unrestricted
nonexistence proof must also cover nonlattice realizations.

A complementary approach studies the distance data rather than the
positions.  Liu \cite{Liu1986} began classifying small configurations
by labeled pair data.  Durst, Hlavacek, Huynh, Miller, and Palsson
\cite{Durst2019} developed this approach using distance matrices and
Cayley--Menger determinants.  Their 2019 version gives three realizable
\emph{strict distance classes} on four points and at least $27$ such
classes on five points.  In the revised manuscript of Chen et al.\
\cite{ChenEtAl2026}, there are $49$ strict five-point classes, exactly
$31$ of which admit planar crescent realizations; the remaining classes
are excluded by exact algebraic calculations.
These counts concern the authors' additional strictness restrictions
on the combinatorial data, not all general-position configurations.
For example, the strict classification excludes a nonrectangular
parallelogram pattern that does admit a general-position realization
\cite{Durst2019,ChenEtAl2026}.  Accordingly, an unrestricted
nonexistence proof must retain such patterns unless a necessary
geometric condition excludes them.

Other settings clarify which hypotheses are essential.
Burt et al.\ \cite{Burt2016} construct $n$-point crescent configurations
in $\R^{n-2}$ for $n\ge4$, with the corresponding hyperplane and
hypersphere general-position conditions.  Thus the multiplicity
sequence alone is not obstructed when the dimension is allowed to
grow.  Fish, King, Miller, Palsson, and Wahlenmayer \cite{Fish2020}
study normed planes, using a strengthened condition that also excludes
four-point distance-equality patterns of equally spaced collinear
points.  They obtain four-point configurations in every norm,
Euclidean examples through size six, and examples through size eight
in the taxicab and Chebyshev norms.
For circular configurations, with distance measured along the shorter
arc, Demaine et al.\ \cite{DemaineEtAl2009} classify the subsets with
this multiplicity sequence, called \emph{Erd\H{o}s-deep} subsets.
Dukes and Gaede \cite{DukesGaede2023} extend this study to families
of modular arithmetic progressions.
These variants concern different geometric constraints and do not
settle the planar Euclidean problem.

The general-position conditions are essential. Clemen, Dumitrescu, and Liu \cite{ClemenDumitrescuLiu2025} construct arbitrarily large planar sets with the required multiplicities by adjoining the center to equally spaced points on a sufficiently short circular arc. These sets are contained in neither a line nor a circle, but when $n\ge5$ they contain four concyclic points and hence are not crescent configurations.

We settle the next size beyond the known planar constructions,
without restricting the coordinates to a lattice or imposing
strictness conditions on the distance pattern.

\begin{theorem}\label{thm:main}
	There is no planar crescent configuration with nine points.
\end{theorem}

Together with the constructions for $2\le n\le8$, the theorem
identifies nine as the first excluded size.  It does not prove
Erd\H{o}s's eventual nonexistence conjecture: deleting points does
not preserve the prescribed multiplicities, so a larger crescent
configuration need not contain a nine-point one.

The proof replaces coordinate search by two finite reductions.
First, the distances of multiplicities one and two occur on just
three pairs.  Deleting one endpoint of each pair, and further points
if necessary, leaves six points whose distances all have original
multiplicity at least three.  We call such a subset a
\emph{high-distance core}; ``high'' refers to multiplicity, not
magnitude.  It determines at most six distances.
We generate all its possible equality patterns, recording which of
the fifteen pairs have equal squared distances, rather than
prescribing a six-point crescent multiplicity sequence.

A finite list of equality patterns need not give finitely many
realizations up to similarity.  We therefore use rank-two Gram
relations and elementary affine identities to obtain necessary
polynomial equations for each pattern.  Rational linear algebra in
the space of homogeneous degree-four polynomials, of dimension at
most $126$, either excludes a pattern or retains it for further
analysis.  A degree-reduction test selects patterns with five or six
distances for which normalization leaves a finite set of candidate
distance tuples.  All other retained patterns are kept in a residual
family.  Compatibility on overlapping six-point subsets then
excludes a nine-point pattern all of whose high-distance cores are
residual.  This forces a core from the selected family without
classifying every continuous family of residual realizations.

Second, multiplication matrices of order at most twelve give
univariate representations of the selected cores.
Once a core is reconstructed, at least two outside points have at
most one low-multiplicity pair joining them to it; these are called
\emph{high-linked points}.  Each is determined either by three
squared distances already present in the core or as the circumcenter
of a core triangle.  After adjoining two such candidates, the last
point has the same property relative to the resulting eight-point
set.  In affine coordinates equipped with their Gram matrix, both
constructions use only field operations.  The extension search
therefore remains in the field generated by the normalized squared
core distances, and exact multiplicity and general-position tests
exclude every candidate.

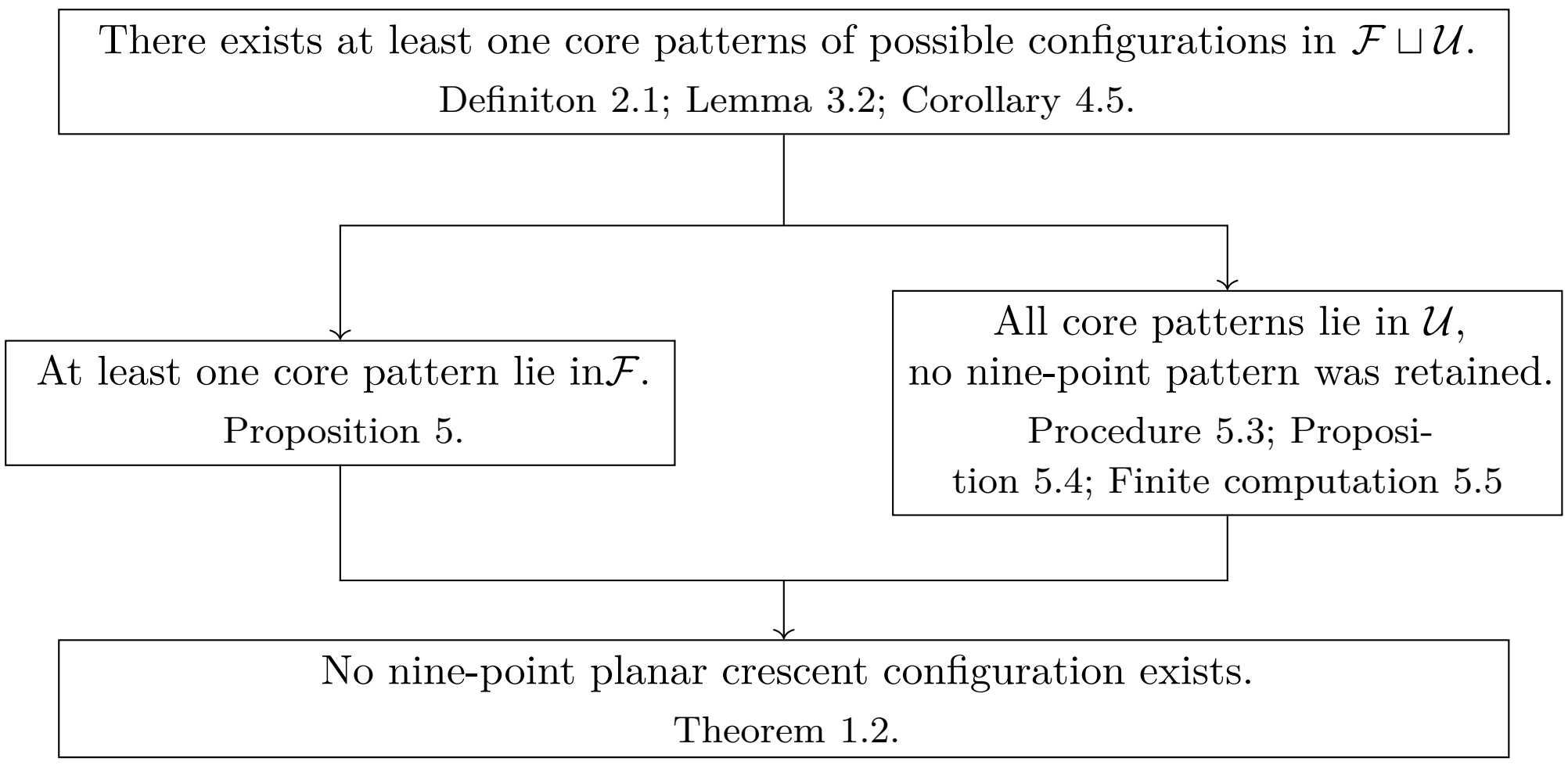
\begin{figure}
\centering
\begin{tikzpicture}[
    x=\linewidth,
    every node/.style={
        draw, align=center, font=\small,
        text width=.80\linewidth,
    }
]

\node (B) at (0,-3.2)    {There exists at least one core patterns of possible configurations in $\mathcal F\sqcup\mathcal U$.\\[2pt]
{\footnotesize
Definiton~\ref{def:pattern};
Lemma~\ref{lem:core};
Corollary~\ref{cor:corepartition}.}};

\node[text width=.36\linewidth] (C) at (-.25,-6) {At least one core pattern lie in$\F$.\\[2pt]
{\footnotesize
Proposition~\ref{sec:coverage}.}};

\node[text width=.36\linewidth] (D) at ( .25,-6) {
All core patterns lie in $\mathcal U$, \\
no nine-point pattern was retained.\\[2pt]
{\footnotesize
Procedure~\ref{proc:overlap};
Proposition~\ref{prop:overlapcomplete};
Finite computation~\ref{comp:overlap}}};

\node (E) at (0,-8.5)     { 
    No nine-point planar crescent configuration exists.\\[2pt]
    {\footnotesize Theorem~\ref{thm:main}.}};

\draw (B.south) -- (0,-4.5);
\draw[->] (0,-4.5) -| (C.north);
\draw[->] (0,-4.5) -| (D.north);
\draw[->] (C.south) -- (-.25,-7.5) -| (E.north);
\draw (D.south) -- (.25,-7.5) -- (0,-7.5);

\end{tikzpicture}
\caption{Dependency structure of the main results.}
\label{structure}
\end{figure}

The distance-matrix viewpoint is classical
\cite{Schoenberg1935,YoungHouseholder1938,Liberti2014}.
The finite-dimensional calculations are related to normal-form
methods \cite{Mourrain1999} and rational univariate representations
\cite{Rouillier1999}; generation up to relabeling belongs to the
general framework of \cite{McKay1998}.
The arguments needed for the present reduction and its completeness
are given explicitly.
Sections~\ref{sec:patterns}--\ref{sec:quartics} generate and classify
the core patterns, Section~\ref{sec:coverage} establishes coverage by
the selected family, and Sections~\ref{sec:fields}--\ref{sec:completion}
give the algebraic reconstruction and exhaustive extension search.
Appendix~\ref{app:implementation} describes the exact evaluation of
these finite procedures. The structure of the proof is presented as Figure~\ref{structure}. 

\section{Distance patterns and elementary obstructions}\label{sec:patterns}

We use squared distances throughout.  For a hypothetical nine-point crescent
configuration $P$, let $\delta_k$ be the squared distance occurring $k$ times.
Thus
\begin{equation}\label{eq:globaldistances}
 \delta_k>0\quad(1\le k\le8),\qquad
 \delta_i\ne\delta_j\quad(i\ne j).
\end{equation}
There is no ordering assumption on these eight numbers.  Pairs carrying $\delta_1$ or $\delta_2$ are called \emph{low}, and pairs carrying one of $\delta_3,\ldots,\delta_8$ are called \emph{high}. These terms refer to occurrence multiplicities. For a finite set
$S\subset\R^2$, write
\[
 \mathcal D(S)=\{|p-q|^2:p,q\in S,\ p\ne q\}.
\]

\subsection{Equality patterns}

We first record which pairs must have equal distances, leaving their
values and the point coordinates to be determined separately.

\begin{definition}[Distance patterns]\label{def:pattern}
	A \emph{distance pattern} on the indices $1,\ldots,n$, with $n\ge2$,
	is a partition of the unordered pairs $\{i,j\}$, $1\le i<j\le n$,
	into nonempty classes.
	Number the classes $1,\ldots,r$ and write $c_{ij}=c_{ji}=\ell$
	when $\{i,j\}$ belongs to class $\ell$.
	The number $r$ of classes is called the \emph{number of labels}.
	Equal labels prescribe equal squared distances, and distinct labels
	prescribe distinct squared distances.  The diagonal has no label.
	
	Two patterns are \emph{equivalent} if one can be obtained from the
	other by permuting the point indices and bijectively renaming the
	labels.
\end{definition}

For algebraic calculations, associate an indeterminate $z_\ell$
with each label and define the \emph{formal squared-distance matrix}
$D=(D_{ij})_{i,j=1}^n$ by
\begin{equation}\label{eq:formalD}
	D_{ii}=0,\qquad D_{ij}=z_{c_{ij}}\quad(i\ne j).
\end{equation}
After assigning values to the variables, we also use $D$ for the
evaluated matrix.

\begin{definition}[Realizations]\label{def:realization}
	A \emph{realization} of a distance pattern is an ordered tuple
	$(p_1,\ldots,p_n)$ of distinct points in $\R^2$ for which pairwise
	distinct positive real values can be assigned to $z_1,\ldots,z_r$
	so that
	\[
	|p_i-p_j|^2=D_{ij}\qquad(1\le i<j\le n).
	\]
	It is \emph{in general position} if its point set is in general
	position as in Definition~\ref{def:crescent}.
\end{definition}

A realization determines the label values uniquely, since every
label occurs.  It has exactly $r$ distinct squared distances,
whose multiplicities are the sizes of the corresponding label
classes.  Its point set is therefore a crescent configuration
precisely when it is in general position and these class sizes
are $1,\ldots,n-1$, each occurring once.
No such multiplicity condition is imposed on a general distance
pattern.  Different realizations of the same pattern need not be
similar, since the pattern prescribes equalities but not distance
ratios.

The \emph{restriction} of a pattern to a subset of its points retains
the labels on pairs within that subset; unused labels are removed and
the remaining ones may be renamed.  For a \emph{partial pattern}, some
pairs are unassigned.  A violation involving only assigned pairs persists
under every completion, which permits the following tests during
generation.

\Needspace{10\baselineskip}
\begin{lemma}\label{lem:local}
The distance pattern of a set in general position has the following
properties.
\begin{enumerate}
 \item At each point, at most three pairs containing it have the same label.
 \item For distinct points $a,b$ and prescribed labels $u,v$, at most two
other points $x$ satisfy $c_{ax}=u$ and $c_{bx}=v$.
 \item For distinct points $a,b$, at most two other points $x$ satisfy
$c_{ax}=c_{bx}$, even when the common label depends on $x$.
\end{enumerate}
\end{lemma}
\begin{proof}
Four equally distant neighbors of a point lie on one circle.  For (ii), the
point $x$ lies on two circles with distinct centers.  Subtracting their
squared-distance equations gives a line, which meets either circle in at
most two points.  In (iii), all such $x$ lie on the perpendicular bisector
of $ab$; three of them would be collinear.
\end{proof}

\begin{lemma}\label{lem:parallelogram}
Let $a,b,c,d$ be four distinct points in general position.  If
\begin{equation}\label{eq:opposite}
 |a-b|=|c-d|,\qquad |a-c|=|b-d|,
\end{equation}
then
\begin{equation}\label{eq:parallelogram}
 a+d=b+c.
\end{equation}
If the third equality $|a-d|=|b-c|$ also holds, the four-point pattern is
impossible in general position.
\end{lemma}
\begin{proof}
After a similarity, take $a=(-1,0)$ and $d=(1,0)$, and write
$b=(u,v)$ and $c=(s,t)$.  Expanding \eqref{eq:opposite} gives
$s=-u$ and $t^2=v^2$.  The choice $t=-v$ gives
\eqref{eq:parallelogram}.  For $t=v$, either $v=0$, giving collinearity,
or all four points lie on the circle centered at
\[
 \left(0,\frac{u^2+v^2-1}{2v}\right).
\]
Thus the latter choice is excluded.  Finally, \eqref{eq:parallelogram}
describes a parallelogram whose diagonals are $ad$ and $bc$.
Equal diagonals make it a rectangle, which is cyclic.
\end{proof}

If $a,b,c$ form an equilateral triangle, their centroid is equidistant
from the three vertices.  The perpendicular bisectors of two sides meet
at a unique point.  Hence any point $d$ equidistant from $a,b,c$ is
that centroid, giving the further affine relation
\begin{equation}\label{eq:centroid}
 a+b+c=3d.
\end{equation}

\subsection{Affine relations and linear distance equations}

An \emph{affine relation} among $p_1,\ldots,p_n$ is an equality
$\sum_i\rho_i p_i=0$ with $\sum_i\rho_i=0$.  For a fixed pattern,
construct a rational matrix $R$ with $n$ columns as follows.  For each
ordered quadruple of distinct indices $(a,b,c,d)$ satisfying
$c_{ab}=c_{cd}$ and $c_{ac}=c_{bd}$, include the coefficient row of
$p_a+p_d-p_b-p_c=0$.  Also include the row of
$p_a+p_b+p_c-3p_d=0$ whenever
$c_{ab}=c_{ac}=c_{bc}$ and $c_{ad}=c_{bd}=c_{cd}$.
These are the relations \eqref{eq:parallelogram} and \eqref{eq:centroid}.
Repeated rows may be omitted; if no relation occurs, $R$ has no rows.
Every realization in general position, with coordinate columns $x,y$,
satisfies
\begin{equation}\label{eq:affinecoordinates}
 R\mathbf1=Rx=Ry=0,
\end{equation}
where $\mathbf1$ is the column of $n$ ones.

\begin{lemma}\label{lem:affine}
A pattern has no realization in general position if the row space of $R$
contains a nonzero vector with at most three nonzero entries.  Moreover,
for every row $\rho=(\rho_i)$ of $R$, every realization in general position
satisfies
\begin{equation}\label{eq:lineardistances}
 \sum_{i=1}^n\rho_i(D_{ij}-D_{i1})=0\qquad(2\le j\le n).
\end{equation}
\end{lemma}
\begin{proof}
Every row-space vector has coordinate sum zero, so a nonzero one has at
least two nonzero entries.  An affine relation with exactly two such
entries forces the corresponding points to coincide.  With exactly three,
at $a,b,c$, the equations
\[
 \rho_a p_a+\rho_b p_b+\rho_c p_c=0,\qquad
 \rho_a+\rho_b+\rho_c=0
\]
give $p_c=p_b-(\rho_a/\rho_c)(p_a-p_b)$, so the three points are collinear.
For the second assertion, evaluate $D$ at a realization in general position
and expand
\[
 \begin{aligned}
 \sum_i\rho_iD_{ij}
 &=\sum_i\rho_i|p_i-p_j|^2\\
 &=\sum_i\rho_i|p_i|^2
   -2p_j\cdot\sum_i\rho_i p_i+|p_j|^2\sum_i\rho_i\\
 &=\sum_i\rho_i|p_i|^2.
 \end{aligned}
\]
The result is independent of $j$, proving \eqref{eq:lineardistances}.
\end{proof}

These conclusions have direct rational linear-algebra tests.
For a triple $I\subseteq\{1,\ldots,n\}$, let $R_{I^c}$ consist of the
columns outside $I$.  Restriction to those columns maps the row space of
$R$ onto the row space of $R_{I^c}$.  It has a nonzero kernel exactly when
\[
 \rank R_{I^c}<\rank R.
\]
Such a kernel vector is a nonzero affine relation supported on $I$.

By \eqref{eq:formalD}, equations \eqref{eq:lineardistances} are homogeneous
linear equations in the label variables.  Collecting their coefficients
gives a rational system
\[
 Az=0,\qquad z=(z_1,\ldots,z_r)^{\mathsf T}.
\]
Let the rows of $N$ be the transposes of a rational basis of $\ker A$.
The transposes of these rows also span the real solution space.  Hence
$Az=0$ forces $z_i=0$ exactly when column $i$ of $N$ is zero, and forces
$z_i=z_j$ exactly when columns $i,j$ agree.  Either condition excludes a
realization in general position, whose label values must be positive and
pairwise distinct.

\section{Six-point cores and exhaustive generation}\label{sec:cores}

The coordinates of a nine-point configuration range over a continuous
space.  We first reduce the search to finitely many equality patterns
on six points, without prescribing their numerical distances.
The reduction uses the fact that the two least frequent distances
occupy only three pairs.

\begin{definition}\label{def:core}
	A \emph{high-distance core} of a nine-point crescent configuration
	$P$ is a six-point subset $S\subset P$ such that
	\begin{equation}\label{eq:core}
		\mathcal D(S)\subseteq\{\delta_3,\ldots,\delta_8\}.
	\end{equation}
\end{definition}

\begin{lemma}\label{lem:core}
	Every nine-point crescent configuration has a high-distance core.
	If the three pairs carrying $\delta_1$ or $\delta_2$ are not pairwise
	disjoint, it has a seven-point subset all of whose squared distances
	belong to $\{\delta_3,\ldots,\delta_8\}$.
\end{lemma}

\begin{proof}
	There is one pair carrying $\delta_1$ and two carrying $\delta_2$.
	Choose one endpoint from each pair and delete the chosen points.
	At most three points are deleted, and none of the three pairs
	remains.  Any six remaining points form a high-distance core.
	
	If two of the pairs share an endpoint, delete that endpoint and
	one endpoint of the third pair.  This removes all three pairs
	using at most two points, so any seven remaining points give
	the second assertion.
\end{proof}

The lemma selects a subset without requiring any knowledge of the
coordinates or distance values.  Its conclusion concerns multiplicities
in $P$: a distance of high multiplicity in $P$ may occur only once
within $S$.  Thus a core need not be a crescent configuration, and
its internal multiplicities are not prescribed.
The seven-point alternative will be used in
Section~\ref{sec:coverage} to compare overlapping six-point cores.

\paragraph{From point sets to a finite list of patterns.}
A six-point set has fifteen unordered pairs.
By \eqref{eq:core}, the pairs of a high-distance core are partitioned
into at most six distance classes.
There are only finitely many such partitions: assigning one of
six labels to each of the fifteen pairs already gives a finite,
though highly redundant, list containing all of them.
We generate these patterns with arbitrary class sizes, without
yet assigning their labels to the global multiplicities
$3,\ldots,8$.

This enumeration does not select points from a grid or approximate
their distances.  It records only which pairs have equal distances.
For each pattern, the distinct positive label values remain to
be determined.  In particular, finiteness of the pattern list
does not assert finiteness of its geometric realizations.

\begin{procedure}[Generation of core patterns]\label{proc:generation}
	Start with the unique two-point pattern, whose pair has label $1$.
	For $n=3,\ldots,6$, extend each retained $(n-1)$-point representative
	by adjoining point $n$ and labeling
	$(1,n),\ldots,(n-1,n)$ in order.
	If the labels currently used are $1,\ldots,t$, try every existing
	label and, when $t<6$, the new label $t+1$.
	
	Reject a partial assignment if its assigned pairs violate
	Lemma~\ref{lem:local} or give the three equal opposite-pair
	comparisons excluded by Lemma~\ref{lem:parallelogram}.
	After completing the new pairs, retain one canonical representative
	of each equivalence class, as specified below.
	After reaching six points, apply the affine and linear-distance
	tests of Lemma~\ref{lem:affine}.
\end{procedure}

All exclusions in this procedure follow from necessary conditions
for general position.  For example, four equal-distance neighbors
force a common circle, regardless of the value assigned to that
distance.  Such a test excludes every general-position realization
of the pattern.  In contrast, the reduction to canonical
representatives removes only different numberings of the same
pattern.

\paragraph{Canonical representatives.}
For any ordering of the points, read the labels in the pair order
\[
(1,2),(1,3),\ldots,(1,n),(2,3),\ldots,(n-1,n),
\]
and rename labels by first occurrence as $1,2,\ldots$.
Compare the resulting sequences lexicographically.

To restrict the point orderings that need to be tested, put
\[
m_\ell=\#\{i<j:c_{ij}=\ell\},\qquad
d_{i\ell}=\#\{j\ne i:c_{ij}=\ell\},
\qquad 1\le\ell\le r,
\]
and associate to each point the sorted list
\begin{equation}\label{eq:pointsignature}
	s_i=\operatorname{sort}\bigl((m_\ell,d_{i\ell})_{\ell=1}^{r}\bigr).
\end{equation}
Sorting is lexicographic and retains repetitions, including
$(m_\ell,0)$ when label $\ell$ does not occur at point $i$.
Group points with equal lists and order the groups by these lists.
Among all point orderings within the groups, take the least
first-occurrence-renamed sequence as the canonical representative.
Retain each resulting sequence once.

\begin{proposition}\label{prop:generation}
	The distance pattern of every six-point set in general position
	with at most six distinct distances is equivalent to a pattern
	retained by Procedure~\ref{proc:generation}.
	In particular, the retained list covers every high-distance core
	of a nine-point crescent configuration.
\end{proposition}

\begin{proof}
	An equivalence preserves the lists \eqref{eq:pointsignature}
	and maps the tested point orderings bijectively to one another.
	Equivalent patterns therefore have the same canonical
	representative.  Conversely, equal representatives give point
	orderings and label renamings that identify the patterns.
	
	For completeness of generation, induct on the number of points.
	The two-point case is immediate.
	Deleting one point from a realizable $n$-point pattern leaves
	a pattern represented at the preceding stage.
	Transfer its points and labels to that representative.
	The pairs joining the deleted point to the others use existing
	labels or new labels, which can be renamed by first occurrence.
	The required assignment is therefore among those tested.
	Every rejection condition is necessary for a realization in
	general position, and canonical reduction removes only equivalent
	copies.  The final affine tests likewise cannot reject a
	general-position realization.
	The assertion about cores follows from Lemma~\ref{lem:core}.
\end{proof}

The output is therefore a finite covering list of possible core
patterns, not a list of verified geometric configurations.
Some retained patterns may have no realization, and a realizable
pattern need not extend to a nine-point crescent configuration.

For each retained pattern with $r\le6$ labels, its fifteen
squared-distance entries in \eqref{eq:formalD} depend only on
$z_1,\ldots,z_r$.  These values may still vary continuously,
even after one is normalized to one.
Section~\ref{sec:quartics} imposes further geometric equations
on them, excluding some patterns and obtaining finite
reconstruction data for others.
The subsequent overlap argument establishes that at least one
high-distance core of any hypothetical nine-point configuration
must belong to this latter family.

\section{Quartic relations and the finite--residual decomposition}
\label{sec:quartics}

Fix a six-point pattern $c$ with $r\le6$ labels retained by Procedure~\ref{proc:generation}. Throughout this section, all realizations under consideration are in general position. We use necessary polynomial relations among the label values to exclude some patterns and to identify others for finite reconstruction after normalization. The retained patterns are then partitioned into $\F$ and $\U$, as specified below.

\subsection{A bounded-degree relation space}

The affine tests have already supplied the linear equations
\eqref{eq:lineardistances}.
To incorporate the requirement that all six points lie in the plane,
form the symmetric $5\times5$ matrix
\begin{equation}\label{eq:Gram}
	L_{ij}=D_{1i}+D_{1j}-D_{ij},
	\qquad 2\le i,j\le6.
\end{equation}
At a realization translated so that $p_1=0$, its entries are
$L_{ij}=2p_i\cdot p_j$.
The evaluated matrix is therefore twice the Gram matrix of five
vectors in $\R^2$ and has rank at most two.
Its $3\times3$ minors give homogeneous cubic equations in the label
variables.  We include the nonprincipal minors.
There are $\binom53=10$ choices of a row triple; symmetry identifies
a minor with the one obtained by interchanging its row and column
triples, leaving $10\cdot11/2=55$ choices.

We place degree-four consequences of the linear and cubic equations
in one coefficient space.  Write $\Q[z]_d$ for the vector space of
homogeneous polynomials of degree $d$ in $z_1,\ldots,z_r$, and set
\begin{equation}\label{eq:quarticdimension}
	V=\Q[z]_4,\qquad
	\dim V=\binom{r+3}{4}\le\binom94=126.
\end{equation}
The dimension is the number of monomials whose nonnegative
exponents sum to four.
The bound $r\le6$, obtained from the core reduction, therefore
limits every coefficient vector to at most $126$ entries.

Let the initial space $W\subseteq V$ be spanned by every cubic
minor multiplied by each variable $z_\ell$, and every linear
relation \eqref{eq:lineardistances} multiplied by each cubic
monomial.  These products all have degree four and vanish at
every general-position realization of $c$.
Their span and membership tests are computed by rational row
reduction in the fixed monomial basis.

The initial span uses polynomial equations but does not yet fully
use the conditions that distances are positive and distinct and
that triples are noncollinear.
These conditions supply the following polynomials that cannot
vanish at a realization:
\begin{equation}\label{eq:nonzero}
	\begin{gathered}
		z_\ell\quad(1\le\ell\le r),\qquad
		z_\ell-z_t\quad(1\le\ell<t\le r),\\
		4D_{ij}D_{ik}-(D_{ij}+D_{ik}-D_{jk})^2
		\quad(1\le i<j<k\le6).
	\end{gathered}
\end{equation}
The first two types are nonzero by positivity and distinctness
of the label values.  At a realization, the last expression is
\[
4D_{ij}D_{ik}-(D_{ij}+D_{ik}-D_{jk})^2
=4\det(p_j-p_i,p_k-p_i)^2>0.
\]
Cancelling any such nonzero factor from a known relation gives a
lower-degree relation.  Its monomial multiples may then add new
relations to $W$ without increasing the working degree.

\begin{lemma}[Division by a known nonzero factor]\label{lem:division}
	Suppose every element of $W\subseteq V$ vanishes at every
	general-position realization of $c$.
	For a factor $f$ in \eqref{eq:nonzero}, let $e=\deg f\in\{1,2\}$
	and set
	\begin{equation}\label{eq:division}
		\begin{aligned}
			C_f&=\{q\in\Q[z]_{4-e}:fq\in W\},\\
			W&\leftarrow
			W+\spanQ\{mq:q\in C_f,\ m\text{ a monomial of degree }e\}.
		\end{aligned}
	\end{equation}
	Every element of the enlarged space still vanishes at every such
	realization.
	Repeated complete passes through the factors terminate when a pass
	leaves $W$ unchanged, with at most $\dim V\le126$ strict enlargements.
\end{lemma}

\begin{proof}
	For $q\in C_f$, the relation $fq=0$ and the condition $f\ne0$
	give $q=0$ at a realization in general position.
	Thus every added polynomial $mq$ vanishes and has degree four.
	Multiplication by the fixed polynomial $f$ is linear in the
	coefficients of $q$, so the condition $fq\in W$ is a rational
	linear system.
	Each strict enlargement increases the dimension of a subspace
	of $V$, proving termination.
\end{proof}

The resulting space consists of necessary relations, not necessarily
all polynomial consequences of the distance equations.
The following tests use these relations to exclude patterns before
attempting to determine their distance values.

\begin{lemma}\label{lem:quarticexclusion}
	A pattern has no realization in general position if $W$ contains
	a degree-four monomial or the determinant of the squared-distance
	matrix of any four of its points.
\end{lemma}

\begin{proof}
	A monomial in positive label values cannot vanish.
	
	For the second assertion, evaluate $D$ at a realization and relabel
	the selected four points as $p_1,\ldots,p_4$, with $p_i=(x_i,y_i)$.
	Let $B$ have rows $(1,x_i,y_i,x_i^2+y_i^2)$.
	Direct multiplication gives
	\begin{equation}\label{eq:circlefactorization}
		(D_{ij})_{i,j=1}^4
		=B\begin{pmatrix}
			0&0&0&1\\
			0&-2&0&0\\
			0&0&-2&0\\
			1&0&0&0
		\end{pmatrix}B^{\mathsf T},
		\qquad
		\det(D_{ij})_{i,j=1}^4=-4(\det B)^2.
	\end{equation}
	If $\det B=0$, a nonzero linear combination of
	$1,x,y,x^2+y^2$ vanishes at all four points.
	The coefficient of $x^2+y^2$ cannot be zero, since no nonzero
	affine linear function vanishes at three noncollinear points.
	Completing squares then gives a circle through the four points;
	its radius is positive because the points are distinct.
	Hence the squared-distance determinant is nonzero in general
	position and cannot belong to $W$.
\end{proof}

\subsection{Certifying finitely many normalized distance tuples}

For a pattern not excluded above, the next task is to determine
whether its distance values can be reduced to a finite calculation.
First remove the common scale: since every label value is positive,
a similarity can normalize any chosen $z_a$ to one.
We then seek relations that replace degree-four expressions by
expressions of degree at most three.

Use $V/W$ as a vector space, identifying two quartics when their
difference belongs to $W$, and write $h=\dim(V/W)$.
Here $h>0$, since $W=V$ would contain a monomial and would already
have excluded the pattern.
For each label $a$, test whether
\begin{equation}\label{eq:finitecriterion}
	V=W+z_a\Q[z]_3.
\end{equation}
This is a rank test: the classes of the quartic monomials divisible
by $z_a$ must span $V/W$.
After setting $z_a=1$, these monomials have degree at most three.
The following lemma extends this reduction to every degree.

\begin{lemma}[Degree reduction]\label{lem:degree}
	Suppose \eqref{eq:finitecriterion} holds.
	Choose degree-four monomials $B_1,\ldots,B_h$, each divisible by
	$z_a$, whose classes form a basis of $V/W$, and put
	$b_i=B_i|_{z_a=1}$.
	Every polynomial in the remaining variables $z_j$, $j\ne a$,
	reduces to a rational linear combination of $b_1,\ldots,b_h$
	using the relations $w|_{z_a=1}=0$, $w\in W$.
	Here a reduction means that the difference is a finite sum of
	polynomial multiples of these relations.
\end{lemma}

\begin{proof}
	A monomial of degree at most four can be multiplied by a power
	of $z_a$ to reach degree four.
	Express its class in $V/W$ in the basis $B_i$, and then set
	$z_a=1$.  This gives a combination of the $b_i$, each of degree
	at most three.
	A monomial of higher degree has a degree-four factor.
	Replace that factor by the preceding expression.
	Every resulting term has smaller total degree, so induction
	completes the reduction.
\end{proof}

The lemma supplies a fixed list of at most $h$ monomials for
calculations at every degree.
Section~\ref{sec:fields} uses the associated multiplication matrices
to obtain explicit finite lists of possible normalized distance
tuples.  Patterns for which this reduction is not available are
not thereby excluded.

\begin{procedure}[Core classification]\label{proc:classification}
	For each six-point pattern passing the affine tests in
	Procedure~\ref{proc:generation}, construct the initial relation
	space $W$ described above.
	Make complete passes through \eqref{eq:division}, applying
	Lemma~\ref{lem:quarticexclusion} to discard a pattern whenever
	either exclusion is obtained.
	Otherwise retain the pattern and its final space $W$ when a
	complete pass leaves $W$ unchanged.
	
	Partition the retained patterns according to their number $r$
	of labels:
	\begin{enumerate}
		\item Assign a pattern to $\F$ if $r\in\{5,6\}$ and
		\eqref{eq:finitecriterion} holds for at least one label $a$.
		Retain the smallest such $a$ and the relation space $W$.
		\item Assign every other retained pattern to $\U$.
		Thus $\U$ includes all retained patterns with at most four
		labels and all five- or six-label patterns for which
		\eqref{eq:finitecriterion} fails for every $a$.
	\end{enumerate}
\end{procedure}

The restriction to five or six labels in $\F$ is used in
Lemma~\ref{lem:twooutside}, which reconstructs two outside points
from such a core.
A retained pattern with at most four labels is assigned to $\U$
even if it satisfies \eqref{eq:finitecriterion}.
Thus $\U$ is the remaining set, not a classification of patterns
with infinitely many realizations.

Procedures~\ref{proc:generation} and \ref{proc:classification} give the counts in Table~\ref{tab:cores}. The family $\U$ consists of three four-label patterns, nine five-label patterns, and $175$ six-label patterns.

\begin{table}[H]
	\centering
	\caption{Core generation and classification, with one representative
		per equivalence class.  The last two rows partition the patterns
		retained by the quartic tests.}
	\label{tab:cores}
	\begin{tabular}{@{}lr@{}}
		\toprule
		Family or stage & Number of patterns\\
		\midrule
		\multicolumn{2}{@{}l}{\emph{Generation with at most six labels}}\\
		\quad Three-point patterns & $3$\\
		\quad Four-point patterns & $22$\\
		\quad Five-point patterns & $932$\\
		\quad Six-point patterns & $625\,472$\\
		\midrule
		\multicolumn{2}{@{}l}{\emph{Filtering and partition of the six-point patterns}}\\
		\quad Retained after the affine tests & $344\,835$\\
		\quad Retained after the quartic tests: $\F\sqcup\U$ & $4\,637$\\
		\qquad Assigned to $\F$ & $4\,450$\\
		\qquad Assigned to $\U$ & $187$\\
		\bottomrule
	\end{tabular}
\end{table}

The generation rows count patterns passing the local tests and the four-point exclusion in Procedure~\ref{proc:generation}, before the six-point affine tests are applied.
The six-point list is then reduced by the affine and quartic tests.
The patterns retained by the quartic tests form exactly the
disjoint union $\F\sqcup\U$; this final partition makes no further
exclusions.  Retention does not assert realizability.
The displayed counts are outputs of the procedures, not conditions
imposed during the search.

\begin{corollary}\label{cor:corepartition}
	The distance pattern of every high-distance core is equivalent
	to a member of $\F\sqcup\U$.
\end{corollary}

\begin{proof}
	A high-distance core is in general position and has at most six
	distinct distances.  Proposition~\ref{prop:generation} therefore
	retains a representative of its pattern after the affine tests.
	Every initial relation in $W$ vanishes at the corresponding
	realization, and Lemma~\ref{lem:division} preserves this property.
	Neither exclusion in Lemma~\ref{lem:quarticexclusion} can apply.
	Its representative is consequently retained and assigned to
	exactly one of $\F$ and $\U$.
\end{proof}

This completes the classification needed from an individual core.
It does not yet show that a core from $\F$ must occur.
Section~\ref{sec:coverage} establishes that assertion by ruling out
a nine-point configuration whose every high-distance core has a
pattern equivalent to a member of $\U$.

\section{A finite core must occur}\label{sec:coverage}

\begin{proposition}[Finite-core coverage]\label{cor:finitecore}
	Every nine-point crescent configuration contains a high-distance core whose distance pattern is equivalent to a member of $\F$.
\end{proposition}

By Corollary~\ref{cor:corepartition}, it suffices to exclude a nine-point crescent configuration in which every high-distance core has a pattern equivalent to a member of $\U$. Under this assumption, the patterns on these overlapping high-distance six-point subsets must be restrictions of a single nine-point pattern.
In the residual search, the labels for $\delta_1$ and $\delta_2$ remain fixed; only high labels may be renamed.
\subsection{Multiplicity capacities}

The \emph{occurrence count} of a label in a partial pattern is the
number of currently assigned unordered pairs carrying that label.
If its final multiplicity is $b$, its current count cannot exceed $b$.
Distinct labels must correspond to distinct final multiplicities.
The following test determines whether such an assignment of
multiplicities is possible.

\begin{lemma}\label{lem:capacity}
  Let $m_1\le\cdots\le m_s$ and $b_1\le\cdots\le b_q$ be positive
  integers, representing current occurrence counts and available final
  multiplicities, respectively.  There exists an injective map
  \[
  \pi:\{1,\ldots,s\}\longrightarrow\{1,\ldots,q\},
  \qquad m_i\le b_{\pi(i)}\quad(1\le i\le s),
  \]
  if and only if
  \begin{equation}\label{eq:capacity}
s\le q
\quad\text{and}\quad
m_i\le b_{q-s+i}\quad(1\le i\le s).
  \end{equation}
\end{lemma}

\begin{proof}
  An injection requires $s\le q$.
  If $m_i>b_{q-s+i}$ for some $i$, the $s-i+1$ counts
  $m_i,\ldots,m_s$ would each require a multiplicity greater than
  $b_{q-s+i}$.  At most $s-i$ available multiplicities are greater
  than this value, so no injection is possible.
  Conversely, under \eqref{eq:capacity}, the assignment
  $\pi(i)=q-s+i$ has the required property.
\end{proof}

In the residual search, the high labels have available multiplicities
$3,4,5,6,7,8$.  Apply Lemma~\ref{lem:capacity} to their current
occurrence counts after each assignment.  The two low labels are
fixed on one and two prescribed pairs, respectively, and cannot be
used on any other pair.  Once all $36$ pairs have been assigned,
require exactly six high labels, with occurrence counts
$3,4,5,6,7,8$.  Before completion, the capacity test is only a
necessary condition: it does not assert that the remaining pairs
can be assigned consistently or that the pattern is realizable.

\subsection{Low-pair arrangements}

Classify the three low pairs up to permutations of the nine point
indices, preserving the pair carrying $\delta_1$.
The three pairs involve distinct endpoints between three and six.
With three endpoints they form a triangle.  With four endpoints
they either share a common endpoint or form a chain of three pairs
through four points.  With five endpoints they consist of two
meeting pairs and a third pair disjoint from both.  With six
endpoints they are pairwise disjoint.  These are the five possible
arrangements before specifying the $\delta_1$ pair.
In a chain, the $\delta_1$ pair can be the middle pair or an end pair.
In the arrangement with two meeting pairs, it can be one of those
two or the disjoint pair.  Within each of the other three
arrangements, all pairs are interchangeable by a permutation of
the endpoints.  Marking the $\delta_1$ pair therefore gives exactly
seven cases.

Table~\ref{tab:low} lists representatives of the seven arrangements,
where $ij$ denotes the unordered pair $\{i,j\}$.
In the first row, every pair among points $1,\ldots,6$ is high;
in each of the other rows, every pair among points $1,\ldots,7$
is high.  Accordingly, Procedure~\ref{proc:overlap} begins with a
pattern from $\U$ on six points in the first case, and with a
seven-point pattern whose every six-point restriction is equivalent
to a member of $\U$ in the other cases.
The last column counts the complete nine-point patterns remaining
after Procedure~\ref{proc:overlap}.  Such a pattern has the prescribed
low pairs and high-label multiplicities $3,\ldots,8$, passes the
local and affine tests specified in that procedure, and has every
six-point restriction containing only high pairs equivalent to a
member of $\U$.  Passing these tests does not by itself establish
geometric realizability.

\begin{table}[H]
  \centering
  \caption{The seven low-pair arrangements and the numbers of complete
nine-point patterns retained by the residual search
in Finite computation~\ref{comp:overlap}.}
  \label{tab:low}
  \begin{tabular}{@{}p{0.50\textwidth}ccc@{}}
\toprule
Low-pair arrangement
    & \shortstack{$\delta_1$\\pair}
    & \shortstack{$\delta_2$\\pairs}
    & \shortstack{Retained\\nine-point\\patterns}\\
\midrule
Three disjoint pairs
    & $47$ & $58,69$ & $0$\\
Triangle
    & $89$ & $18,19$ & $0$\\
Three pairs with a common endpoint
    & $19$ & $29,39$ & $0$\\
Chain; $\delta_1$ on the middle pair
    & $89$ & $18,29$ & $0$\\
Chain; $\delta_1$ on an end pair
    & $18$ & $89,29$ & $0$\\
Two meeting pairs and a disjoint pair;
$\delta_1$ on the disjoint pair
    & $39$ & $18,28$ & $0$\\
Two meeting pairs and a disjoint pair;
$\delta_1$ on a meeting pair
    & $18$ & $28,39$ & $0$\\
\bottomrule
  \end{tabular}
\end{table}
\subsection{Extending residual patterns}

\begin{procedure}[Residual coverage]\label{proc:overlap}
\begin{enumerate}
 \item Adjoin a point to each pattern in $\U$ by the label-assignment rule
of Procedure~\ref{proc:generation}.  Retain a seven-point pattern
only if its restriction to every six-point subset is equivalent to a
member of $\U$.
Apply Lemma~\ref{lem:local} and the exclusion in Lemma~\ref{lem:parallelogram}. For the affine tests, form $R$ using only the parallelogram relations \eqref{eq:parallelogram}, and apply both the row-space and linear-distance tests described after Lemma~\ref{lem:affine}. Apply also Lemma~\ref{lem:capacity} with capacities $3,\ldots,8$, and reduce the retained patterns to canonical representatives.
 \item Fix a row of Table~\ref{tab:low}, setting $m=6$ for the disjoint-pair
row and $m=7$ otherwise.  In the first $m$ positions, place every
pattern in $\U$ if $m=6$, and every seven-point pattern from
\textup{(i)} if $m=7$.  Each placement assigns the points of the
representative to these numbered positions.  All placements are
considered, subject only to the reduction described below.
 \item Adjoin positions $m+1,\ldots,9$ successively.  At each new position,
first prescribe its low pairs from the table.  Assign every remaining
new pair, in increasing order of its old endpoint, an existing high
label or the next unused high label, with at most six high labels in
total.  After an assignment, apply the local and capacity tests and
check each fully assigned six-point restriction containing only high
pairs for membership in $\U$ up to equivalence.  When all pairs to the
new point have been assigned, apply the affine tests from \textup{(i)}.
Reject the current assignment if any of these tests fails, and continue
with the next allowed choice.  At nine points require the low counts
to be $1,2$ and the high counts to be exactly $3,4,5,6,7,8$.
\end{enumerate}
\end{procedure}

\paragraph{Initial placements.}
Fix a low-pair arrangement and place each initial representative on
the first $m$ positions, which contain no low pair.
Group these positions by their prescribed relations to the later
positions: two positions belong to the same group exactly when,
for each later position, both corresponding pairs carry $\delta_1$,
both carry $\delta_2$, or both are high.
Within each group, retain only placements in which the
representative's point numbers increase with the position numbers.
Permuting positions within a group preserves every prescribed low
pair and maps full extensions bijectively to full extensions.
Thus this reduction loses no possible extension.
Finally, rename high labels by first occurrence and remove duplicate
initial sequences, without further permuting their positions.

\paragraph{Generating five new pairs at once.}
In step~\textup{(iii)}, the following lookup can replace the separate assignment of five new high pairs.

For every ordering of every pattern in $\U$, record its restriction to the first five points and the five labels joining them to the sixth. Rename labels by first occurrence, reading the ten old pairs in lexicographic order before the five new pairs in increasing order of their old endpoints. These records form a table indexed by ordered five-point patterns.

When adding a point, choose five old points whose mutual pairs and prescribed pairs to the new point are all high. Order the chosen points by their current indices, and rename the labels on their ten mutual pairs by first occurrence in the same pair order used for the table. Retain the correspondence between these renamed labels and the actual high labels of the current partial pattern.

Look up this ordered five-point pattern. For each matching table entry, use the retained correspondence for labels already present on the ten old pairs. Assign the remaining labels of the entry injectively to high labels absent from those ten pairs. Such labels may already occur elsewhere in the current partial pattern or may be new, provided that at most six high labels are used overall.

Try every such assignment, then fill the remaining pairs and apply the tests of step~\textup{(iii)}. Every extension satisfying the residual-core assumption is represented by a matching table entry, so the absence of a compatible entry excludes the current partial extension under that assumption. If no suitable five-point subset exists, use the individual assignments of step~\textup{(iii)}.

\begin{proposition}\label{prop:overlapcomplete}
If a nine-point crescent configuration has no high-distance core with
pattern equivalent to a member of $\F$,
Procedure~\ref{proc:overlap} retains a complete nine-point pattern.
\end{proposition}
\begin{proof}
By Corollary~\ref{cor:corepartition}, every high-distance core then has a
pattern equivalent to a member of $\U$.  Relabel the low pairs according
to Table~\ref{tab:low}.  In the disjoint case, the initial six points give
an allowed pattern.  Otherwise the first seven points contain only high
pairs, and each six-point restriction is equivalent to a member of $\U$.
Their seven-point pattern is therefore generated in step~\textup{(i)}.
The placement rule in step~\textup{(ii)} keeps a placement related to the
actual one by a permutation preserving the low pairs.

Now adjoin the remaining points in order.  Their actual labels occur
among the choices in step~\textup{(iii)}, including those supplied by the
five-point extension table.  Every completed high-distance restriction
passes the membership test, and every local, affine, and capacity test
is necessary for the actual realization.  Induction therefore retains
its pattern through all nine points, where the prescribed multiplicities
satisfy the final test.
\end{proof}

\begin{computation}\label{comp:overlap}
For the family $\U$ defined by Procedure~\ref{proc:classification},
step~\textup{(i)} of Procedure~\ref{proc:overlap} gives $13$ seven-point
patterns up to equivalence.  In each of the seven arrangements of
Table~\ref{tab:low}, no complete nine-point pattern is retained.
\end{computation}

\begin{proof}[Proof of Proposition~\ref{cor:finitecore}]
	Suppose that a nine-point crescent configuration has no high-distance core whose pattern is equivalent to a member of $\F$. Proposition~\ref{prop:overlapcomplete} then implies that Procedure~\ref{proc:overlap} retains a complete nine-point pattern, contradicting Finite computation~\ref{comp:overlap}.
\end{proof}

\section{Exact distance data for the finite cores}\label{sec:fields}

Fix a pattern in $\F$ and its normalization $z_a=1$. Throughout this section, realizations are understood to be in general position. Choose the monomials $B_i$ in Lemma~\ref{lem:degree} by taking successive independent classes, ordered first by degree after normalization and then lexicographically by their original degree-four exponent vectors. Write $b_i=B_i|_{z_a=1}$.

For a monomial $m$ of degree at most four in the remaining variables, its \emph{coefficient column} is the coordinate vector of the class of $z_a^{4-\deg m}m$ in the basis of $V/W$ formed by the classes of the $B_i$. Extend this definition linearly to polynomials of degree at most four.

Let $e\in\mathbb{Q}^h$ be the coefficient column of the constant
polynomial $1$, and set $M_a=I_h$. For each $j\ne a$, define
$M_j\in\mathbb{Q}^{h\times h}$ by taking its $k$th column to be
the coefficient column of $z_jb_k$, for $1\le k\le h$.
These columns are defined because $\deg b_k\le 3$, so
$\deg(z_jb_k)\le 4$.

For any normalized realization in general position with
squared-distance tuple $(z_1,\ldots,z_r)$, define its evaluation
functional by
\[
    \ell:\mathbb{Q}^h\longrightarrow\mathbb{R},
    \qquad
    \ell(v)=\sum_{k=1}^h v_k b_k(z_1,\ldots,z_r).
\]
The definition of the coefficient columns expresses each
polynomial of degree at most four as a linear combination of
the $b_k$, modulo a relation from $W$ after setting $z_a=1$.
Since every such relation vanishes at the realization, we obtain
\begin{equation}\label{eq:core-evaluation}
    \ell(e)=1,
    \qquad
    \ell(M_jv)=z_j\ell(v)
    \quad
    (1\le j\le r,\ v\in\mathbb{Q}^h).
\end{equation}

The following lemma expresses all normalized squared distances in terms of a single algebraic parameter, giving a finite list of candidate distance tuples.

\begin{lemma}[A univariate certificate]\label{lem:univariate}
	Suppose that the matrices $M_j$ commute and that a rational linear combination $M=\sum_{j=1}^r\alpha_jM_j$ satisfies
	\[
	f(u)=\det(uI-M)\text{ has no repeated root},\qquad
	T=(e,Me,\ldots,M^{h-1}e)\text{ is invertible}.
	\]
	For $1\le j\le r$, define $g_j\in\Q[u]$ by
	\begin{equation}\label{eq:representation}
		g_j(u)=\sum_{k=0}^{h-1}(T^{-1}M_je)_{k+1}u^k.
	\end{equation}
Then $g_j(M)=M_j$ for every $j$. For every normalized realization in general position with squared-distance tuple $(z_1,\ldots,z_r)$, the real number $\lambda=\sum_j\alpha_jz_j$ satisfies
 \begin{equation}\label{eq:univariatetuples}
 	f(\lambda)=0,\qquad z_j=g_j(\lambda)\quad(1\le j\le r).
 \end{equation}
Consequently, the normalized squared-distance tuples of general-position realizations are among the at most $h$ tuples obtained by evaluating $(g_1,\ldots,g_r)$ at the real roots of $f$.
\end{lemma}

\begin{proof}
	The definition of $g_j$ gives
	\[
	g_j(M)e=T(T^{-1}M_je)=M_je.
	\]
	Since both $M_j$ and $g_j(M)$ commute with $M$, we obtain
	\[
	\bigl(g_j(M)-M_j\bigr)M^ke
	=M^k\bigl(g_j(M)-M_j\bigr)e=0
	\qquad(0\le k<h).
	\]
	The vectors $e,Me,\ldots,M^{h-1}e$ form a basis by invertibility of $T$, so $g_j(M)=M_j$.
	
    Now fix a normalized realization with squared-distance tuple
$(z_1,\ldots,z_r)$, and let
$\ell:\mathbb{Q}^h\to\mathbb{R}$ be its evaluation functional
defined above. It satisfies \eqref{eq:core-evaluation}.
	For $\lambda=\sum_j\alpha_jz_j$, these identities give $\ell(Mv)=\lambda\ell(v)$ and hence $\ell(M^ke)=\lambda^k$ for every $k\ge0$. The Cayley--Hamilton theorem gives $f(M)=0$, and therefore
	\[
	f(\lambda)=\ell(f(M)e)=0,\qquad
	g_j(\lambda)=\ell(g_j(M)e)=\ell(M_je)=z_j.
	\]
	Since $f$ has degree $h$ and each real root determines a single candidate tuple, at most $h$ squared-distance tuples can arise from normalized realizations.
\end{proof}

For each pattern in $\F$, choose $M$ by the test in Appendix~\ref{subsec:certificates} and verify the hypotheses of Lemma~\ref{lem:univariate}. Factor $f$ into monic irreducible polynomials over $\Q$. Only factors having a real root can contribute a normalized realization, so discard all other factors.

For each remaining factor $q$, encode the candidate squared distances in
\begin{equation}\label{eq:numberfield}
	K_q=\Q[u]/(q),\qquad z_j=g_j(u)\bmod q.
\end{equation}
Irreducibility of $q$ makes $K_q$ a field. Each element has a unique representative in $\Q[u]$ of degree less than $\deg q$, so equality is tested by comparing these reduced representatives. Sums and products are reduced modulo $q$. Every nonzero representative is coprime to $q$, and the extended Euclidean algorithm therefore computes its inverse modulo $q$. These operations allow the subsequent distance and coordinate calculations to be performed exactly.

Every real root $\theta$ of $q$ defines a field map $K_q\to\R$ by evaluating polynomial representatives at $\theta$. This map is well defined because $q(\theta)=0$, and it is injective because no nonzero rational polynomial of degree less than $\deg q$ can vanish at a root of the irreducible polynomial $q$. We call this map a \emph{real embedding}. Under this embedding, the elements $z_j$ in \eqref{eq:numberfield} become the real values $g_j(\theta)$. Lemma~\ref{lem:univariate} ensures that every normalized realization is represented by one of these factors together with one of its real embeddings.

Injectivity allows zero and equality tests to be performed once for each factor. Discard $q$ if some $z_j$ is zero in $K_q$, or if $z_j=z_k$ for distinct labels $j\ne k$. Such a failure occurs under every real embedding of $K_q$, whereas a realization requires positive, pairwise distinct squared-distance values. Conversely, label values that are nonzero and pairwise distinct in $K_q$ remain so under every real embedding. Their signs may depend on the embedding; positivity and geometric realizability are checked during the reconstruction below.

Finally, suppose that a real root $\lambda$ of a retained factor $q$ represents a normalized realization. Each $g_j(\lambda)$ belongs to $\Q(\lambda)$, while the definition of $\lambda$ gives
\[
\lambda=\sum_j\alpha_jg_j(\lambda).
\]
Since the coefficients of the $g_j$ and the numbers $\alpha_j$ are rational, these two observations give the two inclusions in
\[
\Q(\lambda)=\Q\bigl(g_1(\lambda),\ldots,g_r(\lambda)\bigr).
\]
The image of $K_q$ under evaluation at $\lambda$ is $\Q(\lambda)$. Thus the corresponding embedded field is precisely the field generated by the normalized squared core distances.

\begin{computation}\label{comp:fields}
	Every pattern in $\F$ passes the tests of Lemma~\ref{lem:univariate}, with matrices of order at most $12$. For the choices in Appendix~\ref{subsec:certificates}, the preceding factorization and exclusions retain $4119$ irreducible factors, counted separately for each pattern, before the reconstruction checks below.
\end{computation}

\subsection{Reconstructing a core}\label{subsec:reconstruct}

Fix a pattern in $\F$ and one of its retained factors $q$.
Evaluate the matrix $D$ from \eqref{eq:formalD} at the elements
$z_j\in K_q$ in \eqref{eq:numberfield}.
These are candidate squared distances; we must still determine
whether they have a realization in general position.

In any such realization, take the first point as origin and use
the displacement vectors to the second and third points as a basis.
They are linearly independent by general position.
Writing $p_i$ for the coordinate columns in this basis gives
\begin{equation}\label{eq:affineframe}
	\begin{gathered}
		p_1=\begin{pmatrix}0\\0\end{pmatrix},\qquad
		p_2=\begin{pmatrix}1\\0\end{pmatrix},\qquad
		p_3=\begin{pmatrix}0\\1\end{pmatrix},\\[3pt]
		G=
		\begin{pmatrix}
			D_{12}&(D_{12}+D_{13}-D_{23})/2\\
			(D_{12}+D_{13}-D_{23})/2&D_{13}
		\end{pmatrix}.
	\end{gathered}
\end{equation}
The diagonal entries are the squared lengths of the basis vectors,
and the off-diagonal entry follows by expanding the squared length
of their difference.
Thus $v^{\mathsf T}Gw$ is the inner product of the geometric vectors
with coordinate columns $v,w$.
Write
\[
\|v\|_G^2=v^{\mathsf T}Gv.
\]
The basis need not be orthonormal; recording its metric in $G$
allows the coordinates to be recovered using only operations in
$K_q$.

For this metric to be Euclidean, $G$ must be positive definite
under at least one real embedding.
Discard the factor if $\det G=0$ in $K_q$.
Otherwise test whether some real root of $q$ gives
\begin{equation}\label{eq:positive}
	D_{12}>0,\qquad \det G>0,
\end{equation}
the positive-definiteness criterion for the symmetric matrix $G$.
These inequalities are evaluated at the same root.
Use exact evaluation for a rational root; for any other real root,
refine a rational isolating interval until both signs are determined.
The reduced polynomials representing $D_{12}$ and $\det G$ are
nonzero and therefore cannot vanish at a root of $q$.
Interval refinement consequently terminates.
If no real root satisfies \eqref{eq:positive}, discard the factor.
These are exact univariate operations \cite{BasuPollackRoy2006}.

For a retained factor, recover $p_i=(x_i,y_i)^{\mathsf T}$,
$i=4,5,6$, from its distances to the first three points.
The required equations are
\[
\|p_i\|_G^2=D_{1i},\qquad
\|p_i-p_2\|_G^2=D_{2i},\qquad
\|p_i-p_3\|_G^2=D_{3i}.
\]
Subtracting the first equation from the other two and rearranging
gives
\begin{equation}\label{eq:reconstruct}
	\begin{aligned}
		2D_{12}x_i+(D_{12}+D_{13}-D_{23})y_i
		&=D_{1i}+D_{12}-D_{2i},\\
		(D_{12}+D_{13}-D_{23})x_i+2D_{13}y_i
		&=D_{1i}+D_{13}-D_{3i}.
	\end{aligned}
\end{equation}
The coefficient matrix is $2G$, so each system has a unique
solution in $K_q^2$.
Since these systems use only differences of distance equations,
verify all original distances after reconstruction:
\[
\|p_i-p_j\|_G^2=D_{ij}\qquad(1\le i<j\le6).
\]
Discard the factor if any equality fails.
Finally, require
\begin{equation}\label{eq:generalposition}
	\begin{aligned}
		&\det(p_j-p_i,p_k-p_i)\ne0
		&&(i<j<k),\\
		&\det\begin{pmatrix}
			1&p_i^{\mathsf T}&\|p_i\|_G^2\\
			1&p_j^{\mathsf T}&\|p_j\|_G^2\\
			1&p_k^{\mathsf T}&\|p_k\|_G^2\\
			1&p_l^{\mathsf T}&\|p_l\|_G^2
		\end{pmatrix}\ne0
		&&(i<j<k<l).
	\end{aligned}
\end{equation}
The first determinant has the two vectors as columns and excludes
collinearity.
Under a positive-definite real embedding, changing to orthonormal
coordinates multiplies the second determinant by a nonzero constant
and preserves squared lengths.
Lemma~\ref{lem:quarticexclusion} therefore identifies its vanishing
with concyclicity.

These checks also prove sufficiency.
A positive-definite real matrix $G$ is the Gram matrix of a basis
of $\R^2$.  Interpreting the reconstructed coordinate columns in
that basis realizes the prescribed distances.
They are nonzero, so the points are distinct; positive definiteness
makes their squared distances positive, and
\eqref{eq:generalposition} ensures general position.

\section{Rational completion of a core}\label{sec:completion}

By Proposition~\ref{cor:finitecore}, it remains to determine whether one of the reconstructed six-point sets from $\F$ extends to a nine-point crescent configuration in which it is a high-distance core. We show that every such completion is covered by finite constructions using only field operations in the normalized core-distance field. 

For the geometric arguments below, fix a pattern in $\F$ and a factor $q$ that passes the reconstruction checks of Section~\ref{subsec:reconstruct}. Let $K$ be the image of $K_q$ under a real embedding satisfying \eqref{eq:positive}, and let $S$ be the reconstructed six-point set. In the affine coordinates of \eqref{eq:affineframe}, the points of $S$ lie in $K^2$, and the Gram matrix $G$ has entries in $K$ and is symmetric and positive definite. We write
\[
\|v\|_G^2=v^{\mathsf T}Gv,\qquad
\mathcal D(Q)=\{\|p-q\|_G^2:p,q\in Q,\ p\ne q\}.
\]
Thus $\mathcal D(Q)$ is the squared-distance set in the chosen frame. A squared-distance value is \emph{old} relative to $Q$ if it belongs to $\mathcal D(Q)$.

\begin{lemma}[Three-distance reconstruction]\label{lem:trilateration}
	Let $a,b,c\in K^2$ be noncollinear, and let $\alpha,\beta,\gamma\in K$, not necessarily distinct. At most one point $x\in\R^2$ satisfies
	\[
	\|x-a\|_G^2=\alpha,\qquad
	\|x-b\|_G^2=\beta,\qquad
	\|x-c\|_G^2=\gamma.
	\]
	The only candidate is the point in $K^2$ solving
	\begin{align}
		2(b-a)^{\mathsf T}Gx
		&=b^{\mathsf T}Gb-a^{\mathsf T}Ga+\alpha-\beta,
		\label{eq:trilateration}\\
		2(c-a)^{\mathsf T}Gx
		&=c^{\mathsf T}Gc-a^{\mathsf T}Ga+\alpha-\gamma.\notag
	\end{align}
	It satisfies all three prescribed distances exactly when it satisfies the first squared-distance equation. The circumcenter of $a,b,c$ also lies in $K^2$ and is obtained from the same linear system by setting both distance differences to zero.
\end{lemma}

\begin{proof}
	Subtracting the first squared-distance equation from the other two and rearranging gives \eqref{eq:trilateration}. Its coefficient matrix is the product of the matrix with rows $2(b-a)^{\mathsf T},2(c-a)^{\mathsf T}$ and the matrix $G$. Both are invertible, so the system has a unique solution in $K^2$. If that solution satisfies the first original equation, the two difference equations imply the other two. Setting both distance differences to zero gives the equal-distance equations for the circumcenter.
\end{proof}

This lemma gives a finite construction for a point from an already known set. For a finite set $Q\subset K^2$ in general position, order its points. For each three-point subset, take its vertices in this order and try every ordered triple from $\mathcal D(Q)^3$, allowing repeated values. Solve \eqref{eq:trilateration} and retain the candidate if the first original distance equation holds. Include also the circumcenter of each triangle. Before removing duplicates and points already in $Q$, this uses at most
\begin{equation}\label{eq:recipebound}
	\binom{|Q|}{3}\bigl(|\mathcal D(Q)|^3+1\bigr)
\end{equation}
trials.

For a high-distance core $S$ of a nine-point crescent configuration $P$, a point $x\in P\setminus S$ is \emph{high-linked to $S$} if at least five of its six pairs with points of $S$ are high. Equivalently, at most one of these pairs is low.

\begin{lemma}[Two high-linked points]\label{lem:twooutside}
	Let $S\subset K^2$ be a high-distance core of a nine-point crescent configuration $P$, with $|\mathcal D(S)|\in\{5,6\}$. At least two points of $P\setminus S$ are high-linked to $S$. Every point high-linked to $S$ is obtained either from three squared distances in $\mathcal D(S)$ to a triangle of $S$, or as the circumcenter of a triangle of $S$.
\end{lemma}

\begin{proof}
	There are only three low pairs in $P$. If two of the three points outside $S$ were not high-linked to $S$, each would contribute at least two low pairs joining it to $S$, giving at least four distinct low pairs. Thus at least two outside points are high-linked to $S$.
	
	Fix such a point $x$. If $|\mathcal D(S)|=6$, all high squared distances occur in $S$, so at least five pairs from $x$ to $S$ have old values. Choose three of their endpoints and apply Lemma~\ref{lem:trilateration}.
	If $|\mathcal D(S)|=5$, exactly one high squared distance is absent from $S$. Again, three pairs with old values suffice. Otherwise, among the at least five high pairs from $x$ to $S$, at least three carry the unique high value absent from $S$. Their other endpoints are distinct and noncollinear, so $x$ is the circumcenter of the triangle they form.
\end{proof}

For a core with six squared distances, it suffices to use the two disjoint triangles with indices $(1,2,3)$ and $(4,5,6)$. A point high-linked to $S$ has at most one pair to $S$ with a noncore value, so its three distances to at least one of these triangles are old. Applying the preceding construction to these two triangles gives at most $2(6^3+1)=434$ trials.
For a core with five squared distances, use all $\binom63$ triangles and their circumcenters, giving at most $\binom63(5^3+1)=2520$ trials. In both cases, remove duplicate points and points of $S$ after the prescribed distance checks.

\begin{lemma}[The last point]\label{lem:lastpoint}
	Let $Q\subset K^2$ be an eight-point subset of a nine-point crescent configuration $P$. The remaining point is obtained either from three squared distances in $\mathcal D(Q)$ to a triangle of $Q$, or as the circumcenter of a triangle of $Q$.
\end{lemma}

\begin{proof}
	Write $P\setminus Q=\{x\}$. If some squared distance absent from $\mathcal D(Q)$ has multiplicity at least three in $P$, every pair carrying that distance joins $x$ to $Q$. Choose three of their other endpoints. They are noncollinear and equidistant from $x$, so $x$ is their circumcenter.
	Otherwise, the squared distances absent from $\mathcal D(Q)$ can only be $\delta_1$ and $\delta_2$. Together they account for at most three of the eight pairs from $x$ to $Q$. At least five pairs therefore have old values, and any three of their endpoints reconstruct $x$ by Lemma~\ref{lem:trilateration}.
\end{proof}

\begin{corollary}\label{cor:samefield}
	If a nine-point crescent configuration $P$ has a high-distance core $S\subset K^2$ with five or six squared distances, then every point of $P$ has coordinates in $K^2$, and all squared distances of $P$ lie in $K$. Every such extension is covered by two successive finite constructions: two points high-linked to $S$, followed by the last point relative to the resulting eight-point set.
\end{corollary}

\begin{proof}
	Lemmas~\ref{lem:twooutside} and \ref{lem:trilateration} place two points high-linked to $S$ in $K^2$. All squared distances in the resulting eight-point set therefore lie in $K$. Lemmas~\ref{lem:lastpoint} and \ref{lem:trilateration} then place the last point in $K^2$. The assertion about squared distances follows from the entries of $G$ lying in $K$.
\end{proof}

\paragraph{One search for each irreducible factor.}
The preceding constructions can be carried out directly in the abstract field $K_q$. All candidate coordinates are obtained by field operations. The distance checks and removal of repeated points depend on equality in this field. General position is checked by the nonvanishing of the collinearity determinants and the four-point determinants with rows
\[
(1,p^{\mathsf T},p^{\mathsf T}Gp).
\]
These are the same determinant tests used in Section~\ref{subsec:reconstruct}.

Every real embedding preserves zero and nonzero field elements. It therefore preserves all these tests, as well as equality of squared distances and their occurrence counts. The multiplicity tests below depend only on these counts and integer capacities. For every real embedding satisfying \eqref{eq:positive}, distinct points have positive squared distances because $G$ is positive definite. Consequently, a single search in $K_q$ covers all real embeddings of the factor for which $G$ is positive definite.

\subsection{Multiplicity tests for partial extensions}
\label{subsec:extensiontests}

Let $Q$ consist of a reconstructed six-point set $S$, with $|\mathcal D(S)|\in\{5,6\}$, and one or two candidate outside points. We impose necessary conditions for $Q$ to extend to a nine-point crescent configuration in which $S$ is a high-distance core and the selected outside points are high-linked to $S$.

First discard choices with repeated points or violations of \eqref{eq:generalposition}. For each $d\in\mathcal D(Q)$, let
\[
m_Q(d)=\#\{\{p,q\}\subset Q:p\ne q,\ \|p-q\|_G^2=d\}.
\]
In any completion with the stated properties, the final multiplicities give an injective map $\mu:\mathcal D(Q)\to\{1,\ldots,8\}$ satisfying
\begin{equation}\label{eq:partialcapacity}
	\begin{aligned}
		m_Q(d)&\le\mu(d)
		&&(d\in\mathcal D(Q)),\\
		\mu(d)&\ge3
		&&(d\in\mathcal D(S)),\\
		\#\{p\in S:\mu(\|x-p\|_G^2)\ge3\}&\ge5
		&&(x\in Q\setminus S).
	\end{aligned}
\end{equation}
These conditions respectively bound the current occurrence counts, keep the core distances high, and allow each selected outside point to be high-linked to $S$.

We test the existence of $\mu$ by separating the distances assigned multiplicities $1,2$ from those assigned multiplicities $3,\ldots,8$. Recall from Lemma~\ref{lem:capacity} that sorted occurrence counts $a_1\le\cdots\le a_s$ fit sorted capacities $b_1\le\cdots\le b_t$ exactly when
\[
s\le t,\qquad a_i\le b_{t-s+i}\quad(1\le i\le s).
\]

\paragraph{Six core distances.}
The six core distances must occupy all six high multiplicities. Every distance in $\mathcal D(Q)\setminus\mathcal D(S)$ must therefore have multiplicity $1$ or $2$. Apply Lemma~\ref{lem:capacity} separately to the core occurrence counts with capacities $3,\ldots,8$ and to the noncore occurrence counts with capacities $1,2$.

The last line of \eqref{eq:partialcapacity} becomes
\[
\#\{p\in S:\|x-p\|_G^2\in\mathcal D(S)\}\ge5
\qquad(x\in Q\setminus S).
\]
Thus \eqref{eq:partialcapacity} admits a solution exactly when both capacity tests and this condition hold.

\paragraph{Five core distances.}
Exactly one high distance is absent from the core. For each possible final multiplicity $k\in\{3,\ldots,8\}$ of that distance, first test whether the five core occurrence counts fit the capacities $\{3,\ldots,8\}\setminus\{k\}$. If they do, consider the two possible states of the missing high distance in $Q$:
\begin{enumerate}
	\item It is still absent from $Q$. Then the counts of all observed noncore distances must fit the capacities $1,2$.
	\item It has appeared as some $d_0\in\mathcal D(Q)\setminus\mathcal D(S)$. Try each such $d_0$, require $m_Q(d_0)\le k$, and test whether the counts of all other noncore distances fit the capacities $1,2$.
\end{enumerate}
The observed distances assigned high multiplicities are the core distances, together with $d_0$ in the second case. Require that, for every $x\in Q\setminus S$, at least five pairs from $x$ to $S$ carry one of these high values. The same choice of $k$ and, when present, $d_0$ must work for all selected outside points. Retain $Q$ if at least one choice passes all these tests.

The two cases exhaust the possible assignments in \eqref{eq:partialcapacity}. Within each group, Lemma~\ref{lem:capacity} decides exactly whether the occurrence counts fit the available multiplicities. Passing therefore establishes the existence of a map $\mu$ satisfying \eqref{eq:partialcapacity}, but does not establish the existence of the remaining points.

Existence of $\mu$ is tested afresh for each seven-point and eight-point set. In particular, two candidates that pass separately must still be tested together, since their combined distances and occurrence counts must admit a common assignment.

\begin{procedure}[Complete finite-core extension]\label{proc:extension}
	For every pattern in $\F$ and each of its factors counted in Finite computation~\ref{comp:fields}, do the following.
	\begin{enumerate}
		\item Apply the reconstruction and checks of Section~\ref{subsec:reconstruct}. Discard the factor if any check fails; otherwise let $S$ be the resulting six-point set. Carry out the subsequent extension calculations in $K_q$.
		
		\item Generate candidates for the two high-linked points using the constructions following Lemma~\ref{lem:twooutside}. Remove repetitions and points of $S$. Apply the tests of Section~\ref{subsec:extensiontests} first to each seven-point set obtained by adding one candidate, and then to each eight-point set obtained by adding an unordered pair of distinct retained candidates.
		
		\item For each retained eight-point set $Q$, apply the construction preceding \eqref{eq:recipebound}, using all $\binom83=56$ triangles and the updated squared-distance set $\mathcal D(Q)$. Remove repetitions and points already in $Q$. For each remaining candidate $x$, form $P=Q\cup\{x\}$ and recompute its $36$ squared distances. Retain $P$ exactly when it is in general position, has eight distinct squared-distance values with multiplicities $1,\ldots,8$, and every distance in $\mathcal D(S)$ occurs at least three times.
	\end{enumerate}
\end{procedure}

\begin{proposition}\label{prop:extensioncomplete}
	If a nine-point crescent configuration contains a high-distance core whose pattern is equivalent to a member of $\F$, Procedure~\ref{proc:extension} produces a nine-point crescent configuration.
\end{proposition}

\begin{proof}
	Relabel the core according to its representative in $\F$ and normalize the squared distance selected by \eqref{eq:finitecriterion} to one. Lemma~\ref{lem:univariate} represents its normalized squared-distance tuple by a real root of the characteristic polynomial, and hence by a real embedding of $K_q$ for one of the retained factors $q$.
	
	The reconstruction of Section~\ref{subsec:reconstruct} recovers this core in the frame \eqref{eq:affineframe}. Since the tuple comes from an actual core, all prescribed distances and general-position conditions hold, and the corresponding real embedding satisfies \eqref{eq:positive}. Thus step~\textup{(i)} does not reject its factor.
	Lemma~\ref{lem:twooutside} supplies two points high-linked to the core, both included among the candidates in step~\textup{(ii)}. The seven-point and eight-point subsets formed from these points are in general position, and their actual final multiplicities give maps satisfying \eqref{eq:partialcapacity}. Hence neither the individual candidates nor their pair is rejected.
	Lemma~\ref{lem:lastpoint} supplies the remaining point in step~\textup{(iii)}. The resulting nine-point set passes the full multiplicity, core, and general-position checks.
\end{proof}

\begin{computation}\label{comp:extensions}
	Procedure~\ref{proc:extension}, applied to every pattern and retained factor counted in Finite computation~\ref{comp:fields}, produces no nine-point crescent configuration.
\end{computation}

\begin{proof}[Proof of Theorem~\ref{thm:main}]
	Suppose that a nine-point crescent configuration exists. Proposition~\ref{cor:finitecore} gives a high-distance core whose pattern is equivalent to a member of $\F$. Proposition~\ref{prop:extensioncomplete} would then force Procedure~\ref{proc:extension} to produce a nine-point crescent configuration, contrary to Finite computation~\ref{comp:extensions}.
\end{proof}

\appendix
\section{Exact evaluation of the finite procedures}\label{app:implementation}

\subsection{Shared expansion of the Gram minors}
\label{subsec:coefficientarithmetic}

The initial cubic equations are obtained from a shared expansion.
Assign independent variables to the fifteen pairs of six points and
expand the $55$ Gram minors once.  For each pattern, identify these
variables according to its labels and collect equal monomials.
The resulting integer coefficients are stored for subsequent use.
Since substitution commutes with taking determinants, this avoids
repeated symbolic expansion without changing the equations.

\subsection{Univariate certificates and completion checks}
\label{subsec:certificates}

For each pattern in $\F$, construct $e$ and the matrices $M_j$ as in
Section~\ref{sec:fields} and verify that the matrices commute.
To choose $M$, try $c=0,1,2$ in order: start with $M=0$ and
successively replace $M$ by $cM+M_j$ as $j\ne a$ runs through the
labels in increasing order.  Use the first candidate whose
characteristic polynomial $f(u)=\det(uI-M)$ has no repeated root.
Check that $T$ in Lemma~\ref{lem:univariate} is invertible and verify
the identities $g_j(M)=M_j$.  Finite computation~\ref{comp:fields}
records that these tests succeed for every pattern in $\F$.

Factor $f$ over $\Q$, verifying the factorization and the irreducibility
of its factors.  For each irreducible factor $q$ with a real root,
form $K_q$ and the squared-distance tuple specified in
\eqref{eq:numberfield}.  Discard a factor if any label value is zero
or two label values agree.  The remaining factors give all algebraic
candidates that have not yet been ruled out by positivity or geometry.
Real-root isolation and rational interval refinement determine which
real embeddings satisfy \eqref{eq:positive}; all field arithmetic,
matrix ranks, and equality tests are exact.

The finite calculation is complete only after every generated six-point
pattern has been excluded or assigned to $\F$ or $\U$, all seven
residual searches have been exhausted, and every retained irreducible
factor for every pattern in $\F$ has been checked by
Procedure~\ref{proc:extension}.  A failed certificate leaves its case
unresolved rather than excluding it.  The completed searches retain
no nine-point pattern in the residual case and no nine-point
configuration in the finite-core case.

A complete execution took approximately four hours on a machine with
$16$ GB of RAM configured for $24$ threads.  This is the wall-clock
time of that execution, not a general running-time bound.

\section*{Acknowledgments}

The author thanks Prof.~Eyvindur Palsson for his helpful suggestions.
OpenAI's GPT-6-Astra model assisted in the search for the constructions
presented in this paper. The author takes full responsibility for
the mathematical arguments and results.

\bibliographystyle{plain}
\bibliography{crescent9_v2}
\end{document}